\documentclass{giq9}

\newtheorem*{theorem*}{Theorem}
\newtheorem{prop}{Proposition}
\newtheorem*{prop*}{Proposition}

\newtheorem{rem}{Remark}
\newtheorem{defi}{Definition}
\theoremstyle{definition}
\newtheorem{ex}{Example}

\usepackage{amssymb}

\DeclareMathOperator{\Integer}{\mathbb{Z}}
\DeclareMathOperator{\Complex}{\mathbb{C}}
\DeclareMathOperator{\Real}{\mathbb{R}}
\DeclareMathOperator{\ord}{ord} 
 \DeclareMathOperator{\Ad}{Ad}
\DeclareMathOperator{\CP}{\mathbb{C}P}
\DeclareMathOperator{\HP}{\mathbb{H}P}
\DeclareMathOperator{\Quater}{\mathbb{H}}

\allowdisplaybreaks[1]

\begin{document}

\title[Geometry and Topology  of
Coadjoint Orbits] {GEOMETRY AND TOPOLOGY OF COADJOINT ORBITS OF
SEMISIMPLE LIE GROUPS}
\author{Julia Bernatska$^{\dag \ddag}$, Petro Holod$^{\dag \ddag}$}
\address{$\dag$ Department of Physical and Mathematical Sciences,
National University of 'Kyiv-Mohyla Academy' 2 Skovorody Str., 04070
Kyiv, Ukraine \\ $\ddag$ Bogolyubov Institute for Theoretical
Physics of the National Academy of Sciences of Ukraine 14-b
Metrolohichna Str., 03680, Kyiv, Ukraine }

\begin{abstract}
Orbits of coadjoint representations of classical compact Lie groups
have a lot of applications. They appear in representation theory,
geometrical quantization, theory of magnetism, quantum optics etc.
As geometric objects  the orbits were the subject of much study.
However, they remain hard for calculation and application. We
propose simple solutions for the following problems: an explicit
parameterization of the orbit by means of a generalized
stereographic projection, obtaining a K\"{a}hlerian structure on the
orbit, introducing basis two-forms for the cohomology group of the
orbit.
\end{abstract}

\maketitle

%==============================================================
\section{Introduction}

Orbits of coadjoint representations of semisimple Lie groups are an
extremely interesting subject. These homogeneous spaces are flag
manifolds.  Remarkable, that the coadjoint orbits of compact groups
are K\"{a}hlerian manifolds. In 1950s A.~Borel, R.~Bott,
J.~L.~Koszul, F.~Hirzebruch et al. investigated the coadjoint orbits
as complex homogeneous manifolds. It was proven that each coadjoint
orbit of a compact connected Lie group $\mathrm{G}$ admits a
canonical $\mathrm{G}$-invariant complex structure and the only
(within homotopies) $\mathrm{G}$-invariant K\"{a}hlerian metrics.
Furthermore, the coadjoint orbits can be considered as fibre bundles
whose bases and fibres are coadjoint orbits themselves.

Coadjoint orbits appear in many spheres of theoretical physics, for
instance in representation theory, geometrical quantization, theory
of magnetism, quantum optics. They serve as definitional domains in
problems connected with nonlinear integrable equations (so called
equations of soliton type). Since these equations have a wide
application, the remarkable properties of coadjoint orbits interest
not only mathematicians but also physicists.

It should be pointed out that much of our material is, of course,
not new, but drawn from various areas of the mathematical
literature. The material was collected for solving the physical
problem based on a classical Heisenberg equation with
$\mathrm{SU}(n)$ as a gauge group. The equation describes a behavior
of magnetics with spin $s\,{\geqslant}\, 1$. The paper includes an
investigation of geometrical and topological properties of the
coadjoint orbits. We hope it fulfills a certain need. We would like
to mention that we have added a number of new results (such as an
explicit expression for a stereographic projection in the case of
group $\mathrm{SU}(3)$ and improving the way of its computation, the
idea of obtaining the K\"{a}hlerian potential on an orbit, an
introduction of basis two-forms for the cohomology ring of an
orbit).

The paper is organized as follows. In
section~\ref{s:Coadjoint_orbits} we recall the notion of a coadjoint
orbit, propose a classification of the orbits, and describe the
orbit as a fibre bundle over an orbit with an orbit as a fibre.
Section~\ref{s:Complex_param} is devoted to a generalized
stereographic projection from a Lie algebra onto its coadjoint
orbit, it gives a suitable complex parameterization of the orbit. As
an example, we compute an explicit expression for the stereographic
projection in the case of group $\mathrm{SU}(3)$. In
section~\ref{s:Kahler_structure} we propose a way of obtaining
K\"{a}hlerian structures and  K\"{a}hlerian potentials on the
orbits. Section~\ref{s:topology} concerns a structure of the
cohomology rings of the orbits and finding of G-invariant bases for
the cohomology groups.

%==============================================================
\section{Coadjoint Orbits of Semisimple Lie Groups}
\label{s:Coadjoint_orbits} We start with recalling the notion of a
coadjoint orbit. Let $\mathrm{G}$ be a compact semisimple classical
Lie group, $\mathfrak{g}$ denote the corresponding Lie algebra, and
$\mathfrak{g}^{\ast}$ denote the dual space to $\mathfrak{g}$. Let
$\mathrm{T}$ be the maximal torus of $\mathrm{G}$, and
$\mathfrak{h}$ be the maximal commutative subalgebra (also called a
Cartan subalgebra) of $\mathfrak{g}$. Accordingly,
$\mathfrak{h}^{\ast}$ denotes the dual space to $\mathfrak{h}$.
\begin{defi}
The subset  $\ \mathcal{O}_{\mu} = \{\Ad^{\ast}_g \mu \mid \forall
g\,{\in}\, \mathrm{G} \}$ of $\mathfrak{g}^{\ast}$  is called a
\indd{coadjoint orbit} of $\mathrm{G}$ through~$\mu\in
\mathfrak{g}^{\ast}$.
\end{defi}
In the case of classical Lie groups we can use the standard
representations for adjoint and coadjoint operators:
\begin{equation*}
\Ad_g X = g X g^{-1},\ X\in \mathfrak{g},\qquad \Ad^{\ast}_g \mu =
g^{-1}\mu g,\ \mu\in \mathfrak{g}^{\ast}.
\end{equation*}
Comparing these formulas one can easily see that \emph{a coadjoint
orbit coincides with an adjoint one}.

Define the \emph{stability subgroup} at a point $\mu \in
\mathfrak{g}^{\ast}$ as $\mathrm{G}_{\mu}\,{=}\,\{g\,{\in}\,
\mathrm{G} \mid \Ad_g^{\ast}\mu\,{=}\,\mu \}$. The coadjoint
operator induces a bijective correspondence between an orbit
$\mathcal{O}_{\mu}$ and a coset space $\mathrm{G}_{\mu}\backslash
\mathrm{G}$ (in the sequel, we deal with right coset spaces).

First of all, we classify the coadjoint orbits of an arbitrary
semisimple group $\mathrm{G}$. Obviously, each orbit is drawn from a
unique point, which we call an \emph{initial point} and denote by
$\mu_0$. The following theorem from~\cite{bott79} allows to restrict
the region of search of an initial point.
\begin{theorem*}[R.~Bott]
Each orbit of the coadjoint action of $\mathrm{G}$ intersects
$\mathfrak{h}^{\ast}$ precisely in an orbit of the Weyl group.
\end{theorem*}
In other words, each orbit is assigned to a finite non-empty subset
of $\mathfrak{h}^{\ast}$. For more detail recall the notion of the
Weyl group. Let $N(\mathrm{H})$ be the \emph{normalizer} of a subset
$\mathrm{H}\,{\subset}\, \mathrm{G}$ in $\mathrm{G}$, that is
$N(\mathrm{H})=\{g\,{\in}\, \mathrm{H}\mid g^{-1}\mathrm{H} g =
\mathrm{H}\}$. Let $C(\mathrm{H})$ be the \emph{centralizer} of
$\mathrm{H}$, that is $C(\mathrm{H}) = \{g\,{\in}\, \mathrm{G}\mid
g^{-1}h g \,{=}\, h,\, \forall h\,{\in}\, \mathrm{H}\}$. Obviously,
$C(\mathrm{T})=\mathrm{T}$, where $\mathrm{T}$ is the maximal torus
of $\mathrm{G}$.

\begin{defi}
The \textbf{Weyl group} of $\mathrm{G}$ is the factor-group of
$N(\mathrm{T})$ over $C(\mathrm{T})$:
$$W(\mathrm{G}) = N(\mathrm{T}) / C(\mathrm{T}).$$
\end{defi}
The Weyl group $\mathrm{W}(\mathrm{G})$ acts transitively on
$\mathfrak{h}^{\ast}$. The action of $\mathrm{W}(\mathrm{G})$ is
performed by the coadjoint operator. It is easy to show that
\emph{$\mathrm{W}(\mathrm{G})$ is isomorphic to the finite group
generated by reflections $w_{\alpha}$ across the hyperplanes
orthogonal  to simple roots $\alpha$:
$$w_{\alpha}(\mu) = \mu - 2\tfrac{\langle
\mu, \alpha \rangle}{\langle\alpha,\alpha\rangle}\, \alpha,\qquad
\mu\in \mathfrak{h}^{\ast},$$} where $\langle \cdot,\cdot \rangle$
denotes a bilinear form on $\mathfrak{g}^{\ast}$.

\begin{defi}
The open domain $$C = \{\mu \in \mathfrak{h}^{\ast} \mid \langle
\mu,\alpha \rangle >0, \, \forall \alpha \in \Delta^{+}\}$$ is
called the \indd{positive Weyl chamber}. Here $\Delta^{+}$ denotes
the set of positive roots.

We call the set  $\Gamma_{\alpha} \,{=}\, \{\mu\,{\in}\,
\mathfrak{h}^{\ast} \mid \langle \mu,\alpha \rangle = 0\}$ a
\indd{wall of the Weyl chamber}.
\end{defi}
If we reflect the closure $\overline{C}$ of the positive Weyl
chamber by elements of the Weyl group  we cover
$\mathfrak{h}^{\ast}$ overall:
\begin{equation*}\label{cover_h}
\displaystyle \mathfrak{h}^{\ast} = \bigcup_{w\in
\mathrm{W}(\mathrm{G})} w\cdot \overline{C}.
\end{equation*}

\begin{figure}[h]
\centering \includegraphics[width=10cm]{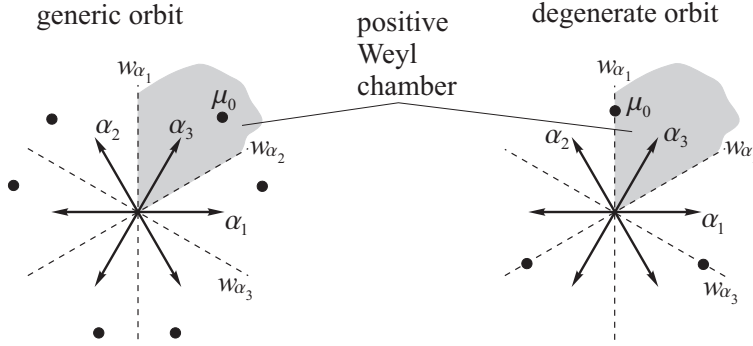} \caption{Root
diagram for $\mathrm{SU}(3)$.\label{f:root_diagr}}
\end{figure}

An orbit of the Weyl group $\mathrm{W}(\mathrm{G})$ is obtained by
the action of $\mathrm{W}(\mathrm{G})$ on a point of $\overline{C}$.
In the case of group $\mathrm{SU}(3)$, two possible types of orbits
of the Weyl group are shown on the root diagram (see
figure~\ref{f:root_diagr}).
Black points denote intersections of a coadjoint orbit with
$\mathfrak{h}^{\ast}$ and form an orbit of
$\mathrm{W(\mathrm{SU}(3))}$. The positive Weyl chamber is filled
with grey color. It has two walls: $\Gamma_{\alpha_1}$ and
$\Gamma_{\alpha_2}$; they are the hyperplanes denoted by
$w_{\alpha_1}$ and $w_{\alpha_2}$. At the left, one can see a
generic case, when an orbit of $\mathrm{W}(\mathrm{\mathrm{SU}(3)})$
has 6 elements. It happens if an initial point lies in the interior
of the positive Weyl chamber. At the right, there is a degenerate
(non-generic) case, when an orbit of
$\mathrm{W}(\mathrm{\mathrm{SU}(3)})$ has 3 elements. It happens if
an initial point belongs to a wall of the positive Weyl chamber.

In the both cases the closed positive Weyl chamber contains a unique
point of an orbit of $\mathrm{W}(\mathrm{G})$. We obtain the
following
\begin{prop}\label{P:classification}
Each orbit $\mathcal{O}$ of $\mathrm{G}$ is uniquely  defined by an
initial point $\mu_0 \in \mathfrak{h}^{\ast}$, which is located in
the closed positive Weyl chamber $\overline{C}$. If $\mu_0$ lies in
the interior of the positive Weyl chamber: $\mu_0\in C$, it gives
rise to a \indd{generic orbit}. If $\mu_0$ belongs to a wall of the
positive Weyl chamber: $\mu_0 \,{\in}\, \Gamma_{\alpha}$,
$\alpha\,{\in}\, \Delta^{+}$, it gives rise to a \indd{degenerate
orbit}.
\end{prop}

As mentioned above, one can define the orbit $\mathcal{O}_{\mu_0}$
through an initial point $\mu_0\,{\in}\, \mathfrak{h}^{\ast}$ by
$\mathcal{O}_{\mu_0} \,{=}\, \mathrm{G}_{\mu_0}\backslash
\mathrm{G}$. Note, that a stability subgroup $\mathrm{G}_{\mu}$ as
$\mu\,{\in}\, \mathfrak{h}^{\ast}$ generically coincides with the
maximal torus $\mathrm{T}$. However, if $\mu$ belongs to a
degenerate orbit, then $\mathrm{G}_{\mu}$ is a lager subgroup of
$\mathrm{G}$ containing $\mathrm{T}$. Therefore, we define a generic
orbit by
$$\mathcal{O}_{\mu_0} = \mathrm{T}\backslash \mathrm{G},$$
and a degenerate one by
$$\mathcal{O}_{\mu_0} = \mathrm{G}_{\mu_0}\backslash \mathrm{G},$$
where $\mathrm{G}_{\mu_0} \neq \mathrm{T}$,
$\mathrm{G}_{\mu_0}\supset \mathrm{T}$.

An important topological property of the coadjoint orbits is the
following. \emph{Almost each orbit can be regarded as a fibre bundle
over an orbit with an orbit as a fibre, except for the maximal
degenerate orbits.} Indeed, if there exists an initial point $\mu_0$
such that $\mathrm{G}_{\mu_0} \,{\supset}\, \mathrm{T}$, one can
form a coset space $\mathrm{T}\backslash \mathrm{G}_{\mu_0}$. Thus,
the orbit $\mathcal{O}_{\mu_0} = \mathrm{T}\backslash \mathrm{G}$ is
a fibre bundle over the base $\mathrm{G}_{\mu_0}\backslash
\mathrm{G}$ with the fibre $\mathrm{T}\backslash
\mathrm{G}_{\mu_0}$:
$$\mathcal{O}_{\mu_0} = \mathcal{E} (\mathrm{G}_{\mu_0}\backslash \mathrm{G}
,\mathrm{T}\backslash \mathrm{G}_{\mu_0},\pi),$$ where $\pi$ denotes
a projection from the orbit onto the base. Moreover,
$\mathrm{G}_{\mu_0}\backslash \mathrm{G}$ and $\mathrm{T}\backslash
\mathrm{G}_{\mu0}$ are coadjoint orbits themselves. We claim this by
\begin{prop}
Suppose $\mathcal{O}_{\mu_0}=\mathrm{G}_{\mu_0}\backslash
\mathrm{G}$  is not the maximal degenerate orbit of $\mathrm{G}$.
Then a subgroup $K$ such that $\mathrm{G}\supset \mathrm{K} \supset
\mathrm{G}_{\mu_0}$ exists, and $\mathcal{O}_{\mu_0}$ is a fibre
bundle over the base $\mathrm{K}\backslash \mathrm{G}$ with the
fibre $\mathrm{G}_{\mu_0}\backslash \mathrm{K}$:
$$\mathcal{O}_{\mu_0} = \mathcal{E} (\mathrm{K}\backslash \mathrm{G},
\mathrm{G}_{\mu_0}\backslash \mathrm{K}, \pi).$$
\end{prop}

We illustrate the proposition by examples.
\begin{ex}
The group $\mathrm{SU}(2)$ has the only type of orbits:
$$\mathcal{O}^{\mathrm{SU}(2)} = \frac{\mathrm{SU}(2)}{\mathrm{U}(1)}
\simeq \CP^1.$$ The group $\mathrm{SU}(3)$ has generic and
degenerate orbits: $$\mathcal{O}^{\mathrm{SU}(3)} =
\frac{\mathrm{SU}(3)}{\mathrm{U}(1)\times \mathrm{U}(1)}, \qquad
\mathcal{O}^{\mathrm{SU}(3)}_{d} =
\frac{\mathrm{SU}(3)}{\mathrm{SU}(2)\times \mathrm{U}(1)} \simeq
\CP^2.$$ Comparing the above coset spaces we see that a generic
orbit $\mathcal{O}^{\mathrm{SU}(3)}$ is a fibre bundle over a
degenerate orbit $\mathcal{O}_d^{\mathrm{SU}(3)}$ with a fibre
$\mathcal{O}^{\mathrm{SU}(2)}$:
$$\mathcal{O}^{\mathrm{SU}(3)}= \mathcal{E}(\mathcal{O}^{\mathrm{SU}(3)}_{d}
, \mathcal{O}^{\mathrm{SU}(2)},\pi) = \mathcal{E} (\CP^2, \CP^1, \pi).$$

The group $\mathrm{SU}(4)$ has several types of degenerate orbits.
There is a list of all possible types of orbits:
\begin{gather*}
\mathcal{O}^{\mathrm{SU}(4)} =
\frac{\mathrm{SU}(4)}{\mathrm{U}(1)\times \mathrm{U}(1)\times
\mathrm{U}(1)}, \qquad \mathcal{O}^{\mathrm{SU}(4)}_{d1} =
\frac{\mathrm{SU}(4)}
{\mathrm{SU}(2)\times \mathrm{U}(1)\times \mathrm{U}(1)}, \\
\mathcal{O}^{\mathrm{SU}(4)}_{d2} =
\frac{\mathrm{SU}(4)}{\mathrm{S}(\mathrm{U}(2)\times
\mathrm{U}(2))}, \qquad \mathcal{O}^{\mathrm{SU}(4)}_{d3} =
\frac{\mathrm{SU}(4)}{\mathrm{SU}(3)\times \mathrm{U}(1)}\simeq
\CP^3.
\end{gather*}
As a result, there exist several representations of a generic orbit
$\mathcal{O}^{\mathrm{SU}(4)}$ as a fibre bundle. For example,
\begin{gather*}
\mathcal{O}^{\mathrm{SU}(4)} = \mathcal{E} (\mathcal{O}^{\mathrm{SU}(4)}_{d3},
\mathcal{O}^{\mathrm{SU}(3)},\pi) = \mathcal{E} (\CP^3,
\mathcal{O}^{\mathrm{SU}(3)},\pi) \\
\mathcal{O}^{\mathrm{SU}(4)} = \mathcal{E} (\mathcal{O}^{\mathrm{SU}(4)}_{d2},
\mathcal{O}^{\mathrm{SU}(2)},\pi) = \mathcal{E} (\mathcal{O}^{\mathrm{SU}(4)}_{d2},
\CP^1,\pi).
\end{gather*}
\end{ex}

\begin{ex}
In the paper we consider compact classical Lie groups. They describe
linear transformations of real, complex, and quaternionic spaces.
Respectively, these groups are $\mathrm{SO}(n)$ over the real field,
$\mathrm{SU}(n)$ over the complex field, and $\mathrm{Sp}(n)$ over
the quaternionic ring.  Here we list the maximal tori of all these
groups, and their representations as fibre bundles.

The maximal torus of $\mathrm{SU}(n)$ is
$\mathrm{T}=\overbrace{\mathrm{U}(1)\times \mathrm{U}(1) \times
\cdots \times \mathrm{U}(1)}^{n-1}$; the generic type of orbits can
be represented as $$\mathcal{O}^{\mathrm{SU}(n)} = \mathcal{E}
(\CP^{n-1}, \mathcal{O}^{\mathrm{SU}(n-1)}, \pi).$$

The maximal torus of $\mathrm{SO}(n)$ as $n=2m$ and $n=2m+1$ has the
following form $\mathrm{T}=\underbrace{\mathrm{SO}(2)\times
\mathrm{SO}(2) \times \cdots \times \mathrm{SO}(2)}_{m}$; the
generic type of orbits can be represented as
\begin{gather*}
\mathcal{O}^{\mathrm{SO}(2m)}= \mathcal{E} (G_{2n;2},
\mathcal{O}^{\mathrm{SO}(2m-2)}, \pi)\\
\mathcal{O}^{\mathrm{SO}(2m+1)}= \mathcal{E} (G_{2n-1;2},
\mathcal{O}^{\mathrm{SO}(2m-1)}, \pi),
\end{gather*}
where $G_{2m; 2}$, $G_{2m-1;2}$ denote real Grassman manifolds.

The maximal torus of $\mathrm{Sp}(n)$ is
$\mathrm{T}=\overbrace{\mathrm{U}(1)\times \mathrm{U}(1) \times
\cdots \times \mathrm{U}(1)}^{n-1}$; the generic type of orbits can
be represented as
$$\mathcal{O}^{\mathrm{Sp}(n)} = \mathcal{E} (\HP^{n-1},
\mathcal{O}^{\mathrm{Sp}(n-1)},\pi),$$ where $\Quater$ denotes the
quaternionic ring.
\end{ex}

%==============================================================
\section{Complex Parameterization of Coadjoint Orbits}\label{s:Complex_param}

In the theory of Lie groups and Lie algebras different ways of
parameterization of coadjoint orbits are available. As the most
prevalent we choose a \indd{generalized stereographic projection}
\cite{skrypnyk}. It is named so since in the case of group
$\mathrm{SU}(2)$ it gives the well-known stereographic projection
onto the complex plane, which is the only orbit of $\mathrm{SU}(2)$.
The generalized stereographic projection is a projection from a dual
space onto a coadjoint orbit parameterized by complex coordinates.

Complex coordinates are introduced by the well-known procedure that
combines Iwasawa and Gauss-Bruhat decompositions. These coordinates
are often called \indd{Bruhat coordinates}~\cite{picken}.

We start with complexifying a group $\mathrm{G}$ in the usual way:
$\mathrm{G}^{\Complex}=\exp\{\mathfrak{g}+\i\mathfrak{g}\}$. A
generic orbit of $\mathrm{G}$ is defined in $\mathrm{G}^{\Complex}$
by Montgomery's diffeomorphism:
\begin{equation}\label{DiffMont}
\mathcal{O} = \mathrm{T}\backslash \mathrm{G} \simeq
\mathrm{P}\backslash \mathrm{G}^{\Complex},
\end{equation}
where  $\mathrm{P}$ denotes the minimal parabolic subgroup of
$\mathrm{G}^{\Complex}$.

Equation \eqref{DiffMont} becomes apparent from the Iwasawa
decomposition $\mathrm{G}^{\Complex} = \mathrm{NAK},$ where
$\mathrm{A}\,{\simeq}\,\exp\{\i\mathfrak{h}\}$ is the real abelian
subgroup of $\mathrm{G}^{\Complex}$, $\mathrm{N}$ is a nilpotent
subgroup of $\mathrm{G}^{\Complex}$, and $\mathrm{K}$ is the maximal
compact subgroup of $\mathrm{G}^{\Complex}$. Since we consider only
compact groups $\mathrm{G}$, $\mathrm{K}$ coincides with
$\mathrm{G}$. Then the Iwasawa decomposition of
$\mathrm{G}^{\Complex}$ has the following form
$$\mathrm{G}^{\Complex} = \mathrm{NAG}.$$

It is easy to express $\mathrm{A}$ and $\mathrm{N}$ in terms of root
vectors. Let $\Delta^{+}$ be the set of positive roots $\alpha$ of
$\mathrm{G}^{\Complex}$. By $X_{\alpha}$, $X_{-\alpha}$, $\alpha\in
\Delta^{+}$, denote positive and negative root vectors,
respectively. By $H_{\alpha}$, $\alpha\,{\in}\, \Delta^{+}$, denote
the corresponding Cartan vectors, which form a basis for the Cartan
subalgebra~$\mathfrak{h}$. According to \cite{hel}, we choose
$X_{\alpha}$ and $X_{-\alpha}$ so that
$X_{\alpha}\,{-}\,X_{-\alpha},\ \i(X_{\alpha}\,{+}\,X_{-\alpha})\in
\mathfrak{g}.$ Then $$\mathrm{N} \simeq
\exp\Bigl\{\textstyle\sum\limits_{\alpha\in \Delta^{+}} n_\alpha
X_{\alpha}\Bigr\},\ n_{\alpha}\in\Complex, \quad \mathrm{A}
\simeq\exp\Bigl\{\textstyle\sum\limits_{\alpha\in \Delta^{+}}
a_\alpha \i H_{\alpha}\Bigr\},\ a_{\alpha}\in \Real.$$

In this notation  $\mathrm{P}=\mathrm{NAT}$. This makes
\eqref{DiffMont} evident.

In the case of a degenerate orbit, we have the following
diffeomorphism:
\begin{equation}\label{DiffMont_deg}
\mathcal{O}_{\mu_0} = \mathrm{G}_{\mu_0} \backslash \mathrm{G}
\simeq \mathrm{P}_{\mu_0} \backslash \mathrm{G}^{\Complex},
\end{equation}
where $\mathrm{G}_{\mu_0}$ is the stability subgroup  and
$\mathrm{P}_{\mu_0}$ is the parabolic subgroup with respect
to~$\mathcal{O}_{\mu_0}$. Then $\mathrm{P}_{\mu_0} =
\mathrm{NAG}_{\mu_0}$, that proves \eqref{DiffMont_deg}.

On the other hand, $\mathrm{G}$ admits a Gauss decomposition (for
the generic type of orbits):
$$\mathrm{G}^{\Complex} =
\mathrm{NT}^{\Complex}\mathrm{Z},$$ where $\mathrm{T}^{\Complex}$ is
the maximal torus of $\mathrm{G}^{\Complex}$, and
$\mathrm{T}^{\Complex}=\mathrm{AT}$ in the above notation;
$\mathrm{N}$ and $\mathrm{Z}\simeq\mathrm{N}^{\ast}$ are nilpotent
subgroups of $\mathrm{G}^{\Complex}$ normalized
by~$\mathrm{T}^{\Complex}$. In terms of the root vectors introduced
above
$$\mathrm{Z} = \exp\Bigl\{\textstyle\sum\limits_{\alpha\in
\Delta^{+}} z_\alpha X_{-\alpha}\Bigr\},\quad z_{\alpha}\in
\Complex.$$

After \cite{hel} we call  $a_{\alpha}$, $n_{\alpha}$,  $z_{\alpha}$
the \emph{canonical coordinates} connected with the root basis
$\{H_{\alpha},\, X_{\alpha},\, X_{-\alpha} \mid \alpha\in
\Delta^{+}\}$. These are coordinates in the group $\mathrm{G}$.

A comparison of the Gauss and Iwasawa decompositions implies that
the orbit $\mathcal{O}$ is diffeomorphic to the subgroup manifold
$\mathrm{Z}$:
\begin{equation}\label{GSProj}
\mathcal{O}\simeq \frac{\mathrm{NAG}}{\mathrm{NAT}} \simeq
\frac{\mathrm{NATZ}}{\mathrm{NAT}}\simeq\mathrm{Z}.
\end{equation}
Diffeomorphism \eqref{GSProj} asserts that one can parameterize the
orbit $\mathcal{O}$ in terms of the complex coordinates
$\{z_{\alpha}$, $\alpha\in \Delta^{+}\}$ that are canonical
coordinates in $\mathrm{Z}$.

However, a Gauss decomposition is local. Therefore, we use a
Gauss-Bruhat decomposition instead:
\begin{gather*}
\mathrm{G}^{\Complex} = \bigcap_{w\in \mathrm{W}(\mathrm{G})}
\mathrm{P} \mathrm{Z}w.
\end{gather*}
It gives a system of local charts on the orbit:
\begin{equation}\label{Maps}
\mathcal{O} = \mathrm{P} \backslash \mathrm{G}^{\Complex} =
\bigcap_{w\in \mathrm{W}(\mathrm{G})} \mathrm{Z}w.
\end{equation}

In the case of a degenerate orbit $\mathcal{O}_{\mu_0}$,
$\mathrm{T}$ is to be replaced by $\mathrm{G}_{\mu_0}$, and
$\mathrm{P}$  by $\mathrm{P}_{\mu_0}$. It is sufficient to take the
intersection over $w\in \mathrm{W}(\mathrm{G}_{\mu_0})\backslash
\mathrm{W}(\mathrm{G})$ in \eqref{Maps}. Furthermore, in this case,
$\mathrm{Z}$ has a less number of coordinates.

\begin{prop}
Each orbit $\mathcal{O}$ of a compact semisimple Lie group
$\mathrm{G}$ is locally parameterized in terms of the canonical
coordinates $\{z_{\alpha}$, $\alpha \,{\in}\, \Delta^{+}\}$ in a
nilpotent subgroup~$\mathrm{Z}$ of $\mathrm{G}^{\Complex}$ according
to \eqref{Maps}.
\end{prop}

Now we apply the above scheme to compact classical Lie groups,
namely $\mathrm{SO}(n)$, $\mathrm{SU}(n)$, $\mathrm{Sp}(n)$. The
scheme consists of several steps. First we parameterize the
subgroups $\mathrm{N}$, $\mathrm{A}$, and the group $\mathrm{G}$ in
terms of $\{z_{\alpha}$, $\alpha\in \Delta^{+}\}$. Secondly, we
choose an initial point $\mu_0$ in the positive closed Weyl chamber
$\overline{C}$ and generate an orbit $\mathcal{O}_{\mu_0}$ by the
dressing formula
$$\mu = g^{-1}\mu_0 g,\quad g\in \mathrm{G}.$$ That gives a
parameterization on one of the charts covering the orbit. Finally,
we extend the parameterization to all other charts by the action of
elements of the Weyl group of~$\mathrm{G}$. We consider the scheme
in detail.

\textbf{Step 1}. Being a finite group, each classical Lie group has
a matrix representation. Let $\hat{a}$ be the matrix representing an
element $a$. An Iwasawa decomposition of $\hat{z}\in \mathrm{Z}$ has
the following form:
\begin{equation}\label{decompZ}
\hat{z} = \hat{n}\hat{a}\hat{k},\qquad \hat{n}\in \mathrm{N},\quad
\hat{a}\in \mathrm{A},\quad \hat{k}\in \mathrm{G}.
\end{equation}
One has to solve \eqref{decompZ} in terms of the complex coordinates
$z_{\alpha}$ that appear as entries of the matrix $\hat{z}$. The
following transformation of \eqref{decompZ} makes the computation
easier
$$\hat{z}\hat{z}^{\ast} =
\hat{n}\hat{a}\hat{k}\hat{k}^{\ast}\hat{a}^{\ast}\hat{n}^{\ast} =
\hat{n}\hat{a}^2\hat{n}^{\ast},$$ where $\hat{k}^{\ast}$ denotes the
hermitian conjugate of $\hat{k}$. Indeed, $\hat{k}\hat{k}^{\ast} =
e$ for all of the mentioned groups. This is evident, if one
considers the conjugation over the complex field in the case of
$\mathrm{SU}(n)$, and over the quaternionic ring in the case of
$\mathrm{Sp}(n)$. If $\hat{k}\in \mathrm{SO}(n)$ one has
$\hat{k}^{\ast}=\hat{k}^T$, and the equality $\hat{k}\hat{k}^{\ast}
= e$ is obvious. Moreover, it can easily be checked  that
$\hat{a}\hat{a}^{\ast}=\hat{a}^2$. When $\hat{n}$ and $\hat{a}$ are
parameterized in terms of  $\{z_{\alpha}\}$, the matrix $\hat{k}(z)$
is computed by the formula
$$\hat{k}(z) = \hat{a}^{-1}(z)\hat{n}^{-1}(z)\hat{z}.$$

Here we obtain complex parameterizations of $\mathrm{N}$,
$\mathrm{A}$, $\mathrm{G}$ for all classical compact groups of small
dimensions.

\begin{ex}
In the case of group $\mathrm{SU}(n)$, the corresponding
complexified group is $\mathrm{SL}(n,\Complex)$. The subgroup
$\mathrm{N}$ consists of complex upper triangular matrices with ones
on the diagonal, the subgroup $\mathrm{Z}$ consists of complex low
triangular matrices with ones on the diagonal, the subgroup
$\mathrm{A}$ contains real diagonal matrices
$\hat{a}=\diag(r_1,\,r_2,\,\dots,\, r_n)$ such that $\prod_{i=1}^n
r_i = 1$.

Decomposition \eqref{decompZ} for a generic orbit
$\mathcal{O}^{\mathrm{SU}(3)}$ gets the form
\begin{gather*}
\small \begin{pmatrix} 1&0&0 \\ z_1&1&0 \\
z_3&z_2&1 \end{pmatrix} = \begin{pmatrix} 1&n_1&n_3 \\
0&1&n_2 \\ 0&0&1\end{pmatrix} \begin{pmatrix} \frac{1}{r_1}&0&0 \\ 0&\frac{r_1}{r_2}&0 \\
0&0&r_2 \end{pmatrix}\hat{u}, \qquad \hat{u} \in \mathrm{SU}(3),
\end{gather*}
whence it follows
\begin{gather*}
r^2_1 = 1+|z_1|^2+|z_3-z_1z_2|^2,\quad r^2_2 = 1+|z_2|^2 + |z_3|^2\\
n_1 = \frac{1}{r_1^2}(\bar{z}_1(1+|z_2|^2)-z_2\bar{z}_3),\quad n_2=
\frac{1}{r_2^2}(\bar{z}_2+z_1 \bar{z}_3),\quad n_3 =
\frac{\bar{z}_3}{r_2^2}.
\end{gather*}
The dressing matrix $\hat{u}$ is $$\hat{u}=\begin{pmatrix}
\frac{1}{r_1} &
-\frac{\bar{z}_1}{r_1} & -\frac{\bar{z}_3-\bar{z}_1 \bar{z}_2}{r_1}\\
\frac{z_1(1+|z_2|^2)-z_3\bar{z}_2}{r_1r_2} &
\frac{1+|z_3|^2-z_1z_2\bar{z}_3}{r_1 r_2} & -\frac{\bar{z}_2 +
z_1\bar{z}_3}{r_1 r_2} \\ \frac{z_3}{r_2} & \frac{z_2}{r_2} &
\frac{1}{r_2}\end{pmatrix}.$$

The case of a degenerate orbit $\mathcal{O}_d^{\mathrm{SU}(3)}$ is
derived from the above by assigning $z_1\,{=}\,0$, or $z_2\,{=}\,0$.
\end{ex}

\begin{ex}
In the case of group $\mathrm{Sp}(n)$, the complexified group is
$\mathrm{Sp}(n,\Complex)$. The both groups describe
 linear transformations of the quaternionic vector
space $\Quater^n$. Therefore, it is suitable to operate with
quaternions instead of complex numbers. Each quaternion $q$ is
determined by two complex numbers $z_1$, $z_2$ as
$q=z_1+z_2\mathbf{j}$. The quaternionic conjugate of $q$ is
$\bar{q}=\bar{z}_1-\mathbf{j}\bar{z}_2$, where $\bar{z}_1$,
$\bar{z}_2$ are the complex conjugates of $z_1$, $z_2$. Several
useful relations are available: $$\mathbf{j}z =
\bar{z}\mathbf{j},\qquad \overline{z + w} = \bar{z}+\bar{w},\qquad
\overline{z\cdot \vphantom{+} w} = \bar{w}\cdot  \bar{z},$$ where
$z,\, w\in \Complex$.

The subgroups $\mathrm{N}$, $\mathrm{Z}$ have the same
representatives as in the case of group $\mathrm{SU}(n)$, but over
the quaternionic ring. The subgroup $\mathrm{A}$ consists of real
diagonal matrices with the same property as in the case of
$\mathrm{SU}(n)$.

We start with the simplest group $\mathrm{Sp}(2)$. Suppose $v,\,
q\in \Quater$ such that $v=n_1+n_2\mathbf{j}$,
$q=z_1+z_2\mathbf{j}$, where $n_1,\, n_2,\,z_1,\, z_2\in\Complex$.
Decomposition \eqref{decompZ} for an orbit
$\mathcal{O}^{\mathrm{Sp}(2)}$ gets the following form
\begin{gather*}
\small \begin{pmatrix} 1&0 \\ q&1  \end{pmatrix} = \begin{pmatrix} 1&v \\
0&1 \end{pmatrix} \begin{pmatrix} \frac{1}{r}&0 \\ 0&r
\end{pmatrix}\hat{p}, \qquad \hat{p} \in \mathrm{Sp}(2),
\end{gather*}
whence it follows $r^2 = 1+|q|^2$, $v=\bar{q} / r^2$, or in terms of
complex coordinates:
\begin{gather*}
r^2 = |z_1|^2 + |z_2|^2,\quad n_1 = \frac{\bar{z}_1}{r^2},\quad n_2=
-\frac{z_2}{r^2}.
\end{gather*}

The dressing matrix $\hat{p}$ is
$$\hat{p}=\frac{1}{\sqrt{|z_1|^2 + |z_2|^2}} \begin{pmatrix}
 1 & -\bar{z}_1+\mathbf{j} \bar{z}_2 \\
z_1 + z_2\mathbf{j} & 1 \end{pmatrix}.$$

\bigskip

In the case of group $\mathrm{Sp}(3)$, we perform all computations
in terms of quaternions. Suppose $q_1=z_1+z_2\mathbf{j}$,
$q_2=z_3+z_4\mathbf{j}$, $q_3=z_5+z_6\mathbf{j}$,
$v_1=n_1+n_2\mathbf{j}$, $v_2=n_3+n_4\mathbf{j}$,
$v_3=n_5+n_6\mathbf{j}$. Then, for a generic orbit
$\mathcal{O}^{\mathrm{Sp}(3)}$, one obtains
\begin{gather*}
\small \begin{pmatrix} 1&0&0 \\ q_1&1&0 \\
q_3&q_2&1 \end{pmatrix} = \begin{pmatrix} 1&v_1&v_3 \\
0&1&v_2 \\ 0&0&1\end{pmatrix} \begin{pmatrix} \frac{1}{r_1}&0&0
\\ 0&\frac{r_1}{r_2}&0 \\
0&0&r_2 \end{pmatrix}\hat{p}, \qquad \hat{p} \in \mathrm{Sp}(3),
\end{gather*}
whence it follows
\begin{gather*}
r^2_1 = 1+|q_1|^2+ |q_3-q_2 q_1|^2,\quad r^2_2 = 1+|q_2|^2 + |q_3|^2 \\
v_1 = \frac{1}{r_1^2}(\bar{q}_1(1+|q_2|^2)-\bar{q}_3 q_2),\quad v_2=
\frac{1}{r_2^2}(\bar{q}_2 + q_1 \bar{q}_3),\quad v_3 =
\frac{\bar{q}_3}{r_2^2}.
\end{gather*}

The dressing matrix $\hat{p}$ is
$$\hat{p}=\begin{pmatrix}
\frac{1}{r_1} &
-\frac{\bar{q}_1}{r_1} & -\frac{\bar{q}_3-\bar{q}_1 \bar{q}_2}{r_1}\\
\frac{q_1(1+|q_2|^2)-\bar{q}_2 q_3}{r_1r_2} &
\frac{1+|q_3|^2-q_1\bar{q}_3q_2}{r_1 r_2} & -\frac{\bar{q}_2 +
q_1\bar{q}_3}{r_1 r_2} \\ \frac{q_3}{r_2} & \frac{q_2}{r_2} &
\frac{1}{r_2}\end{pmatrix}.$$

The case of $\mathrm{Sp}(n)$ in terms of quaternions is very similar
to the case of $\mathrm{SU}(n)$. The only warning is that the
multiplication of quaternions is not commutative.
\end{ex}

\begin{ex}
In the case of group $\mathrm{SO}(n)$, the corresponding
complexified group is $\mathrm{SO}(n,\Complex)$. Representatives of
the subgroups $\mathrm{N}$ and $\mathrm{Z}$ have not so clear
structure as for groups $\mathrm{SU}(n)$ and $\mathrm{Sp}(n)$. The
real abelian subgroup $\mathrm{A}$ consists of block-diagonal
matrices $\hat{a} = \diag(A_1,\,A_2,\,\dots,\, A_m)$ in the case of
group $\mathrm{SO}(2m)$, and $\hat{a} = \diag(A_1,\,A_2,\,\dots,\,
A_m,\,1)$ in the case of group $\mathrm{SO}(2m+1)$. Here $$A_i= \begin{pmatrix} \cosh a_i & -\i\sinh a_i \\
\i\sinh a_i & \cosh a_i \end{pmatrix}.$$

Consider the group $\mathrm{SO}(3)$.  The only type of orbits is
$\mathcal{O}^{\mathrm{SO}(3)}\,{=}\,\mathrm{SO}(2)\backslash
\mathrm{SO}(3)$. In this case decomposition \eqref{decompZ}  gets
the form
\begin{gather*}
\begin{pmatrix} 1-\frac{z^2}{2} & -\frac{\i z^2}{2} & -z
\\ -\frac{\i z^2}{2}& 1+\frac{z^2}{2} & -\i z \\ z & \i z & 1
\end{pmatrix} = \begin{pmatrix} 1-\frac{n^2}{2}&\frac{\i n^2}{2} & n \\
\frac{\i n^2}{2}& 1+\frac{n^2}{2} & -\i n \\ -n & \i n & 1
\end{pmatrix} \begin{pmatrix} \cosh a & -\i\sinh a & 0 \\ \i\sinh a & \cosh a & 0\\
0&0&1 \end{pmatrix} \hat{o},
\end{gather*}
where $\hat{o}\,{\in}\,\mathrm{SO}(3)$, and $a$, $n$, $z$ are
canonical coordinates in the group. One easily computes the
following
\begin{gather*}
e^{a} = 1+|z|^2,\qquad n=\frac{\bar{z}}{1+|z|^2}.
\end{gather*}
The dressing matrix $\hat{o}$ is
$$\hat{o}=\begin{pmatrix}
\frac{2-z^2-\bar{z}^2}{2(1+|z|^2)} &
\frac{\i(\bar{z}^2-z^2)}{2(1+|z|^2)} & -\frac{z+\bar{z}}{1+|z|^2}\\
\frac{\i(\bar{z}^2-z^2)}{2(1+|z|^2)} &
\frac{2+z^2+\bar{z}^2}{2(1+|z|^2)}  & -\frac{\i(z-z\bar{z})}{1+|z|^2} \\
\frac{z+z\bar{z}}{1+|z|^2} & \frac{\i(z-\bar{z})}{1+|z|^2} &
\frac{1-|z|^2}{1+|z|^2}\end{pmatrix}.$$
\end{ex}

We return to the scheme.

\textbf{Step 2}. Suppose we have some parameterization of the dual
space $\mathfrak{g}^{\ast}$ to the algebra~$\mathfrak{g}$ of a group
$\mathrm{G}$. We call these parameters group coordinates. In order
to parameterize an orbit of $\mathrm{G}$ we find expressions for the
group coordinates in terms of the complex coordinates
$\{z_{\alpha}$, $\alpha\in\Delta^{+}\}$. Continue the example of
group $\mathrm{SU}(3)$.

Let $\lambda_a$, $a=1..\,8$, be Gell-Mann matrices, then $Y_a =
-\frac{\i}{2}\lambda_a$, $a=1..\,8$, form a basis
for~$\mathfrak{g}^{\ast}$. Define a bilinear form on
$\mathfrak{g}^{\ast}$ as $\langle A,B \rangle = -2\Tr AB$. Each
basis element $Y_a$ is assigned to a group coordinate: $\mu_a =
\langle \hat{\mu}, Y_a \rangle,$ where
$$\hat{\mu}=-\frac{\i}{2} \small\begin{pmatrix}
\mu_3+\frac{1}{\sqrt{3}}\mu_8 & \mu_1-\i\mu_2& \mu_4-\i\mu_5 \\
\mu_1+\i\mu_2 & -\mu_3+\frac{1}{\sqrt{3}}\mu_8& \mu_6-\i\mu_7 \\
\mu_4+\i\mu_5&\mu_6+\i\mu_7& -\frac{2}{\sqrt{3}}\mu_8
\end{pmatrix}.$$

A coadjoint orbit is generated by the dressing formula:
$$\hat{\mu} = \hat{u}^{\ast}\hat{\mu}_0\hat{u},\qquad \hat{\mu}_{0}\in\mathfrak{h}^{\ast},$$ where
$\hat{\mu}_0$ is an initial point. As shown in
section~\ref{s:Coadjoint_orbits}, each orbit is uniquely defined by
a point of the closed positive Weyl chamber. Let simple roots of
$\mathfrak{su}(3)$  be as follows: $\hat{\alpha}_1 = \diag
(\i,-\i,0)$ and $\hat{\alpha}_2 = \diag(0,\i,-\i)$. The closed
positive Weyl chamber is the set of points $\hat{\mu}_0$ such that
\begin{equation}\label{init_point}
\hat{\mu}_0 = - \frac{\i}{3}\,\xi\begin{pmatrix} 2&0&0 \\
0&-1 & 0 \\ 0&0&-1\end{pmatrix}
-\frac{\i}{3}\,\eta\begin{pmatrix} 1&0&0 \\
0&1 & 0 \\ 0&0&-2\end{pmatrix},\qquad \xi,\, \eta > 0.
\end{equation}
Obviously, walls of the Weyl chamber are obtained by assigning
$\xi=0$ or $\eta=0$. In this notation $\Gamma_{\alpha_1} =
\{-\frac{\i}{3}\eta\diag(1,1,-2) \mid \eta >0\}$,
 $\Gamma_{\alpha_2} = \{-\frac{\i}{3}\xi\diag(2,-1,-1)
\mid \xi>0\}$. The chosen representation of an initial point
$\hat{\mu}_0$ is the most suitable for the further computation.

According to Proposition~\ref{P:classification} we get a generic
orbit if $\eta\neq 0$ and $\xi\neq 0$. If $\xi$ or $\eta$ vanishes,
we get a degenerate one. A generic orbit is parameterized by three
complex coordinates $z_1$, $z_2$, $z_3$. If $\xi$ vanishes, one has
to  assign $z_1\,{=}\,0$. If $\eta$ vanishes, then  $z_2=0$. We
consider the degenerate orbit through the following point
\begin{equation*}
\hat{\mu}_0 =  -\frac{\i}{3}\,\eta \begin{pmatrix} 1&0&0 \\
0&1 &0 \\ 0&0&-2\end{pmatrix}.
\end{equation*}

One can attach some physical meaning to nonzero entries of the
initial point $\hat{\mu}_0$ because of its diagonal form. For in
quantum mechanics  diagonal matrices represent observable variables.
Suppose $\hat{\mu}_0$ is the value of $\hat{\mu}$ at the infinity:
$\hat{\mu}_0=\hat{\mu}(\infty)$. The diagonal entries are expressed
in terms of the group coordinates $\mu_3$ and $\mu_8$; we fix their
values at the infinity: $\mu_3(\infty)=m$, $\mu_8(\infty)=q$. Then
$$\eta=-\tfrac{1}{2}\left(m-\sqrt{3}q\right),\qquad \xi= m.$$ Suppose
the group $\mathrm{SU(3)}$ describes a magnetic with spin 1. Then
$m$ serves as a projection of magnetic moment (magnetization) of the
magnetic, and $q$ serves as a projection of quadrupole moment.

The dressing procedure gives the following explicit expression for
the generalized stereographic projection onto a generic orbit of
$\mathrm{SU}(3)$:
\begin{equation}\label{SU3_StereoP}
\begin{split}
&\mu_1 = -\frac{\eta}{r_2^2} (\bar{z}_2z_3+z_2\bar{z}_3)
-\frac{\xi}{r_1^2}(z_1+\bar{z}_1)\\
&\mu_2 = \phantom{-}\frac{\i\eta}{r_2^2} (\bar{z}_2z_3-z_2\bar{z}_3)
+\frac{\i\xi}{r_1^2}(z_1-\bar{z}_1)\\
&\mu_3 = \phantom{-}\frac{\eta}{r_2^2}(|z_2|^2-|z_3|^2)
+\frac{\xi}{r_1^2}(1-|z_1|^2) \\
&\mu_4 =  -\frac{\eta}{r_2^2}  (z_3+\bar{z}_3)
-\frac{\xi}{r_1^2}(z_3-z_1z_2+\bar{z}_3-\bar{z}_1\bar{z}_2)\\
&\mu_5 =  \phantom{-} \frac{\i\eta}{r_2^2}  (z_3-\bar{z}_3)
+\frac{\i\xi}{r_1^2}\bigl(z_3-z_1z_2-(\bar{z}_3-\bar{z}_1\bar{z}_2)\bigr) \\
&\mu_6 =  -\frac{\eta}{r_2^2}  (z_2+\bar{z}_2)
+\frac{\xi}{r_1^2}\bigl(\bar{z}_1(z_3-z_1z_2)+z_1(\bar{z}_3-\bar{z}_1\bar{z}_2)\bigr)\\
&\mu_7= \phantom{-} \frac{\i\eta}{r_2^2}  (z_2-\bar{z}_2)
-\frac{\i\xi}{r_1^2}\bigl(\bar{z}_1(z_3-z_1z_2)-z_1(\bar{z}_3-\bar{z}_1\bar{z}_2)\bigr)\\
&\sqrt{3}\mu_8 =\frac{\eta}{r_2^2}(2-|z_2|^2-|z_3|^2)
+\frac{\xi}{r_1^2} (1+|z_1|^2-2|z_3-z_1z_2|^2),
\end{split}
\end{equation}
where $$r_1^2 = 1+|z_1|^2+|z_3-z_1z_2|^2,\quad r_2^2
=1+|z_2|^2+|z_3|^2.$$ Obviously, all expressions can be divided into
two parts: with the coefficients $\eta$ and $\xi$. These parts
correspond to the basis matrices in \eqref{init_point}.

For the stereographic projection  onto a degenerate orbit through
$\hat{\mu}_0$ chosen above one has to assign $\xi=0$, $z_1=0$ in
\eqref{SU3_StereoP}.

\textbf{Step 3}. Parameterization \eqref{SU3_StereoP} is available
on the coordinate chart containing the point $(z_1=0,\, z_2=0,\,
z_3=0)$. By the action of elements of the Weyl group  one obtains
parameterizations on all other charts. The Weyl group is generated
by reflections  across the hyperplanes orthogonal to simple roots.
In the case of group $\mathrm{SU}(3)$, these reflections are
represented by the following matrices
$$\hat{w}_1 = \begin{pmatrix} 0&1&0 \\ 1&0&0 \\
0&0&-1 \end{pmatrix},\qquad \hat{w}_2 = \begin{pmatrix} -1&0&0 \\
0&0&1 \\ 0&1&0 \end{pmatrix}.$$ The action of $\hat{w}_1$ transforms
the chart with coordinates \eqref{SU3_StereoP} onto another one by
the following change of coordinates: $$(z_1,z_2,z_3) \mapsto (z_1',
z_2',z'_3),\qquad z_1'=\frac{1}{z_1},\ z_2'=-z_3,\ z'_3=-z_2.$$ This
chart contains the point $(z_1=\infty,\, z_2=0,\, z_3=0)$. The
action of $\hat{w}_2$ transforms coordinates \eqref{SU3_StereoP} by
the following change of coordinates: $$(z_1,z_2,z_3) \mapsto (z_1',
z_2',z'_3),\qquad z_1'=-(z_3-z_1 z_2),\ z_2'=\frac{1}{z_2},\
z'_3=-\frac{z_3}{z_2}.$$ The latter chart contains the point
$(z_1=0,\, z_2=\infty,\, z_3=0)$.

Evidently, the other elements of $\mathrm{W}(\mathrm{SU}(3))$ are
$\hat{e}$, $\hat{w}_1\hat{w}_2$, $\hat{w}_2\hat{w}_1$,
$\hat{w}_1\hat{w}_2\hat{w}_1$. The corresponding changes of
coordinates are obtained by sequential actions of the two described
above.

%==============================================================
\section{K\"{a}hlerian Structure on Coadjoint Orbits}
\label{s:Kahler_structure}

The perfect property of coadjoint orbits of compact semisimple Lie
groups is the following. Each orbit is simultaneously a Riemannian
manifold and a symplectic one. A Riemanian metrics and the matched
symplectic form together are called a \indd{K\"{a}hlerian
structure}. A.~Borel \cite{borel} proved the following
\begin{prop} Suppose $\mathrm{G}$ is a semisimple compact Lie
group. Then each orbit of~$\mathrm{G}$ admits a complex analytic
K\"{a}hlerian structure invariant under the group $\mathrm{G}$.
\end{prop}

It means that each orbit possesses a hermitian Riemannian metrics,
the K\"{a}hlerian metrics $ds^2$, and the corresponding closed
two-form, the \indd{K\"{a}hlerian form} $\omega$:
$$ds^2 = \sum_{\alpha,\beta} g_{\alpha\bar{\beta}} dz_{\alpha} d\bar{z}_{\beta},\qquad
\omega = \sum_{\alpha,\beta} \i g_{\alpha\bar{\beta}}\, dz_{\alpha}
\wedge d\bar{z}_{\beta}.$$

The $\mathrm{G}$-invariance of a K\"{a}hlerian structure means
invariance under the action of~$\mathrm{G}$. Here we consider the
action of a group as right multiplication. A K\"{a}hlerian structure
is determined by a \indd{K\"{a}hlerian potential} $\Phi$ according
to the formula
$$g_{\alpha\bar{\beta}} = \frac{\partial^2 \Phi}{\partial z_{\alpha} \partial
\bar{z}_{\beta}},\qquad\qquad \omega_{\alpha\bar{\beta}} = \i
g_{\alpha\bar{\beta}}.$$

The objective of this section is to obtain an expression for a
K\"{a}hlerian structure on a coadjoint orbit. Evidently, for this
purpose it is sufficient to find a K\"{a}hlerian potential, which
simultaneously gives the K\"{a}hlerian metrics and the K\"{a}hlerian
form.

On the other hand, one has the following
\begin{prop}[see \cite{besse}]
If $\mathrm{G}$ is a compact semisimple Lie group, the
Kirillov-Kostant-Souriau two-form coincides with a
$\mathrm{G}$-invariant K\"{a}hlerian form.
\end{prop}
While we deal with compact semisimple classical Lie groups, we can
use a Kirillov-Kostant-Souriau differential form as a K\"{a}hlerian
form.

Define a bilinear form on $\mathfrak{g}$ as follows:
$$\langle X,Y \rangle = \Tr XY,\qquad  X,\, Y \in \mathfrak{g}.$$
In the case of classical Lie groups, the bilinear form is
proportional to the standard Killing form on $\mathfrak{g}$.

Define  a vector field $\widetilde{X}$ on a coadjoint orbit
$\mathcal{O}$ by
$$\widetilde{X} f(\mu) = \frac{d}{d\tau} \left. f(\Ad^{\ast}_{\exp(\tau X)} \mu )\right|_{\tau=0},
\quad f\in C^{\infty}(\mathcal{O}).$$ One can introduce an
$\Ad$-invariant closed two-form on $\mathcal{O}$ by the formula
\begin{equation}\label{KKSform}
\omega(\widetilde{X},\widetilde{Y}) = \langle \mu,
[X,Y]\rangle,\qquad X,Y\in \mathfrak{g},\quad
\mu\in\mathfrak{g}^{\ast}.
\end{equation}
This two-form is called a \indd{Kirillov-Kostant-Souriau form}.

The straightforward way of obtaining a K\"{a}hlerian form is to
solve equations \eqref{KKSform}. Unfortunately, it becomes extremely
complicate in dimensions greater than 3. This way is developed by
R.~F.~Picken in~\cite{picken}. He computes K\"{a}hlerian forms on
flag manifolds  via $\mathrm{G}$-invariant one-forms in terms of
Bruhat coordinates.

We return to the idea of finding a K\"{a}hlerian potential instead
of a K\"{a}hlerian form. In general, each $\mathrm{G}$-covariant
real function on an orbit serves as a K\"{a}hlerian potential. It
turns out, that each orbit has a unique $\mathrm{G}$-covariant real
function, which we call a \indd{K\"{a}hlerian potential on the
orbit}.

The same idea is used by D.~V.~Alekseevsky and A.~M.~Perelomov
in~\cite{aleks}. In order to find potentials for all closed
two-forms on orbits of group $\mathrm{GL}(n)$, they consider the
real positive functions built by means of principal minors of
$\hat{z}\hat{z}^{\ast}\,{\in}\, \mathrm{GL}(n)$, and select the
functions that are $\mathrm{G}$-covariant. Here we develop the idea
of D.~V.~Alekseevsky and A.~M.~Perelomov, because this way allows to
avoid complicate computations.

Below we prove that a K\"{a}hlerian potential is determined by a
one-dimensional irreducible representation of the real abelian
subgroup $\mathrm{A}$ of $\mathrm{G}^{\Complex}$. We use the
group-theoretical approach in our proof.

Each orbit $\mathcal{O}= \mathrm{P} \backslash
\mathrm{G}^{\Complex}$ is a holomorphic manifold, which admits the
construction of a line bundle. Let $\{U_{k}\}$ be its atlas. An
arbitrary $g_{\Complex} \in \mathrm{G}^{\Complex}$ has a
decomposition
\begin{equation}\label{sect}
g_{\Complex} = h_{k}(x) s_{k}(x),\qquad x\in U_{k},
\end{equation}
where $s_{k}: U_{k} \to \mathrm{G}^{\Complex}$ is a local section of
$\mathcal{O}$.  If $U_{k} \cap U_{j} \neq \emptyset$, then there
exists a map $s_{kj}=s_{k} \circ s_{j}^{-1}$, which is $s_{kj}:
U_{k} \cap U_{j} \to \mathrm{P}$. A one-dimensional representation
of the parabolic subgroup $\mathrm{P}$ of $\mathrm{G}^{\Complex}$
gives  a $\mathrm{G}$-covariant function on an orbit.

Recall, that $\mathrm{P} = \mathrm{N}\mathrm{A}\mathrm{T}$ in the
case of a generic orbit. In the case of a degenerate orbit, one has
$\mathrm{P} = \mathrm{N}\mathrm{A}\mathrm{G}_{\mu_0}$, where
$\mathrm{G}_{\mu_0}$ is the stability subgroup at an initial point
$\mu_0\in \mathfrak{h}^{\ast}$ giving rise to the orbit. A
one-dimensional irreducible representation is trivial on any
nilpotent group. This means that the representation of $\mathrm{P}$
coincides with the representation of the maximal torus
$\mathrm{T}^{\Complex}\,{=}\,\mathrm{AT}$ of
$\mathrm{G}^{\Complex}$. Moreover, we are interested in real
representations because a K\"{a}hlerian potential is a real
function. Consequently, the required representation is determined
only by $\mathrm{A}$.

Now we build a one-dimensional irreducible representation of
$\mathrm{T}^{\Complex}$. Obviously, $\mathrm{T}^{\Complex}$ is
isomorphic to a direct product of $l$ samples of the multiplicative
group $\Complex^{\ast}\,{=}\,\Complex\backslash\{0\}$, where
$l\,{=}\,\dim \mathrm{T}$. Let the following set of complex numbers
$(d_1$, $d_2,\, \dots,\, d_l)$ be an image of $\hat{d}\in
\mathrm{T}^{\Complex}$ under the isomorphism. It is clear that the
set of real numbers $(r_1,\,r_2,\, \dots,\, r_l)$, where
$r_i=|d_i|$, $i\,{=}\,1..\, l$, is an image of $\hat{a}\,{\in}\,
\mathrm{A}$ under the isomorphism. In terms of complex coordinates
$z=\{z_{\alpha}\mid\alpha\,{\in}\,\Delta^{+}\}$, which are canonical
coordinates in $\mathrm{Z}$, an Iwasawa decomposition of any
$\hat{z}\in \mathrm{Z}$ gets the form
\begin{equation}\label{expIw}
\hat{z} = \hat{n}(z)\hat{a}(z)\hat{k}(z).
\end{equation}
Here $\hat{k}(z)$ represents a point of an orbit in terms of the
complex coordinates $\{z_{\alpha}\}$; $\hat{n}(z)$ and $\hat{a}(z)$
denote matrices $\hat{n}$ and $\hat{a}$ in terms of
$\{z_{\alpha}\}$. After the action of an element $g\in\mathrm{G}$ on
$z$ we perform a Gauss-Bruhat decomposition:
\begin{equation}\label{expGB}
\hat{z} \hat{g} =\hat{n}_B(zg) \hat{d}(zg) \hat{z}_g,\qquad
\hat{n}_B(zg)\in \mathrm{N}.
\end{equation}

From the Iwasawa decomposition of $\hat{z}_g$ we have
$$\hat{a}(z_g) = \hat{n}^{-1}(z_g)\hat{z}_g \hat{k}^{-1}(z_g).$$
Using \eqref{expIw} and \eqref{expGB} we get
\begin{equation}\label{f1}
\hat{a}(z_g) = \hat{n}^{-1}(z_g)\hat{d}^{-1}(zg)\,
\hat{n}^{-1}_B(zg)\hat{n}(z)\hat{a}(z)\hat{k}(z) \hat{g}
\hat{k}^{-1}(z_g).
\end{equation}
In order to gather nilpotent elements together we recall that the
maximal torus $\mathrm{T}^{\Complex}$ is the normalizer of
$\mathrm{N}$, that gives the following equality $$\hat{d}^{-1}(zg)
\hat{n}_B^{-1}(zg) \hat{n}(z)\hat{d}(zg)= \hat{n}(z,g), \qquad
\hat{n}(z,g)\in \mathrm{N}.$$ Substituting
$\hat{n}_B(zg)\hat{d}^{-1}(zg)$ for $\hat{d}^{-1}(zg)
\hat{n}_B^{-1}(zg) \hat{n}(z)$ in \eqref{f1} we obtain
\begin{equation*}
\hat{a}(z_g) =
\hat{n}^{-1}(z_g)\hat{n}(z,g)\hat{d}^{-1}(zg)\hat{a}(z)\hat{k}(z)
\hat{g} \hat{k}^{-1}(z_g).
\end{equation*}
To cancel the element $\hat{k}(z) \hat{g}
\hat{k}^{-1}(z_g)=\hat{g}'\in \mathrm{G}$ we take the following
product
\begin{equation}\label{f2}
\hat{a}^2(z_g) = \hat{a}(z_g)\hat{a}^{\ast}(z_g) = \hat{n}
\hat{d}^{-1}(z_g)\hat{a}^2(z)
\hat{d}^{\ast}{}^{-1}(z_g)\hat{n}^{\ast},
\end{equation}
where $\hat{n}$ denotes $\hat{n}^{-1}(z_g)\hat{n}(z,g) \in
\mathrm{N}$.

Now we construct a one-dimensional real representation of
\eqref{f2}. Let $\chi^{\xi}(\hat{a})$ denote a representation of
$\hat{a}$ with real weights $\xi=(\xi_1,\,\xi_2,\,\dots,\, \xi_l)$.
A one-dimensional real representation of
$\hat{a}\,{\in}\,\mathrm{A}$ has the following form
$\chi^{\xi}(\hat{a}) \,{=}\, r_1^{\xi_1}r_2^{\xi_2}\cdots
r_l^{\xi_l}$, and a one-dimensional real representation of
$\hat{d}\,{\in}\,\mathrm{T}^{\Complex}$  has the form
$\chi^{\xi}(\hat{d}) = d_1^{\xi_1}d_2^{\xi_2}\cdots d_l^{\xi_l}$.
Therefore, the representation of $\hat{a}^2(z_g)$ gets the form
$$\chi^{2\xi}\bigl(\hat{a}(z_g)\bigr)  =
\chi^{\xi}\bigl(\hat{d}(zg)\bigr)
\overline{\chi^{\xi}\bigl(\hat{d}(zg)\bigr)}
\chi^{2\xi}\bigl(\hat{a}(z)\bigr).$$ Whence it is seen that
$\chi^{2\xi}\bigl(\hat{a}(z)\bigr)$ is transformed by a cocycle
$\chi^{\xi}\bigl(\hat{d}(zg)\bigr)$ defined on $\mathrm{G} \times
\mathcal{O}$. It means that the function
\begin{equation}\label{Kpotential}
\ln \chi^{2\xi}\bigl(\hat{a}(z)\bigr) = \xi_1 \ln r_1^2(z) + \xi_2
\ln r_2^2(z) + \dots + \xi_l \ln r_l^2(z)
\end{equation}
is $\mathrm{G}$-covariant, and serves as a K\"{a}hlerian potential
on $\mathcal{O}$. Moreover, each function $\ln r_i^2(z)$,
$i\,{=}\,1..\, l$, is a K\"{a}hlerian potential itself.

Remarkably, that each coadjoint orbit has a unique K\"{a}hlerian
potential of the form \eqref{Kpotential}, where the weights
$\xi=(\xi_1,\,\xi_2,\,\dots,\, \xi_l)$ are determined by an initial
point of the orbit. We have proven the following
\begin{prop}
Suppose $\mathrm{A}$ is the real abelian subgroup of
$\mathrm{G}^{\Complex}$, $\hat{a}\in \mathrm{A}$, and
$\chi^{\xi}(\hat{a})$ is a one-dimensional representation of
$\hat{a}$ with real weights $\xi=(\xi_1,\,\xi_2,\,\dots,\, \xi_l)$.
Then K\"{a}hlerian potentials on coadjoint orbits of $\mathrm{G}$
have the form $\ln \chi^{2\xi}(\hat{a})$, moreover each orbit has
the K\"{a}hlerian potential with a unique $\xi$.
\end{prop}

\begin{rem}
In the case of integer weights $\xi=(\xi_1,\,\xi_2,\,\dots,\,
\xi_l)$, the line bundle over each coadjoint orbit of $\mathrm{G}$
is holomorphic. This idea is derived from the Borel-Weyl theory
based on \cite{borel_weyl}.
\end{rem}

Consider several examples.
\begin{ex}
In the case of group $\mathrm{SU}(n)$, a representative of the real
abelian subgroup $\mathrm{A}$ has the
 form of a diagonal matrix with $\det \hat{a} = 1$, that is
 $$\hat{a} = \diag(1/r_1,\, r_1/ r_2,\,\dots,\, r_{n-2} / r_{n-1},\,
r_{n-1}),$$ and $\dim \mathrm{T}\,{=}\,n\,{-}\,1$. Let $(r_1,\,
r_2,\,\dots,\, r_{n-1})$ be an image of $\hat{a}$ under an
isomorphism from $\mathrm{T}^{\Complex}$ onto
$(\Complex^{\ast})^{n-1}$. Then $\chi^{\xi}(\hat{a}) =
r_1^{\xi_1}r_2^{\xi_2}\cdots r_{n-1}^{\xi_{n-1}}$, where
$\xi_i\,{\in}\,\Real$, $i\,{=}\,1..\,(n\,{-}\,1)$, whence
K\"{a}hlerian potentials have the following form: $$\Phi = \xi_1 \ln
r_1^2 + \xi_2 \ln r_2^2 + \cdots + \xi_{n-1} \ln r_{n-1}^2.$$

For instance, K\"{a}hlerian potentials on orbits of $\mathrm{SU}(3)$
are
\begin{gather*}
\Phi =  \xi \ln r_1^2 + \eta \ln r_2^2 \\
r_1^2 = 1+|z_1|^2+|z_3-z_1z_2|^2,\qquad r^2_2 = 1+|z_2|^2 + |z_3|^2.
\end{gather*}
This expression completely accords with the straightforward solution of
\eqref{KKSform}, which gives the following:
\begin{gather*}
\Phi = \langle \hat{\mu}_0,\hat{\alpha}_1 \rangle\Phi_1 +
\langle\hat{\mu}_0,\hat{\alpha}_2 \rangle \Phi_2\\
\Phi_1=\ln (1+|z_1|^2+|z_3-z_1 z_2|^2), \qquad \Phi_2=
\ln(1+|z_2|^2+|z_3|^2),
\end{gather*}
here $\hat{\mu}_0$ is an initial point of an orbit,
$\hat{\alpha}_1$, $\hat{\alpha}_2$ are the simple roots of
$\mathfrak{su}(3)$. In the case of a degenerate orbit, one has to
assign $z_1=0$ or $z_2=0$.
\end{ex}

\begin{ex}
In the case of groups $\mathrm{SO}(n)$, $n=2m$ and $n=2m+1$, a representative of the subgroup $\mathrm{A}$ has the
form of a block-diagonal matrix, namely
$$\begin{array}{l} \phantom{or\ }\hat{a} = \diag(A_1,\, A_2,\,\dots,\, A_m)
\\ \text{or\ } \hat{a} = \diag(A_1,\, A_2,\,\dots,\, A_m,1) \end{array},\quad A_i =
\begin{pmatrix} \cosh a_i & -\i\sinh a_i \\
\i\sinh a_i & \cosh a_i \end{pmatrix},\ i=1..\, m.$$ Here $\{a_i\}$
are canonical coordinates in the maximal torus $\mathrm{T}$, and
$\dim \mathrm{T}=m$. Let $(e^{a_1},\, e^{a_2},\,\dots,\, e^{a_m})$
be an image of $\hat{a}$ under an isomorphism from
$\mathrm{T}^{\Complex}$ onto $(\Complex^{\ast})^{m}$. Then
$\chi^{\xi}(\hat{a}) = e^{\xi_1 a_1}e^{\xi_2 a_2}\cdots e^{\xi_m
a_m},$ whence it follows $$\Phi = 2\xi_1 a_1 + 2\xi_2 a_2 + \cdots
+2\xi_m a_m.$$

K\"{a}hlerian potentials on coadjoint orbits of $\mathrm{SO}(4)$
computed by \eqref{KKSform} have the form
\begin{gather*}
\Phi = \langle \hat{\mu}_0,\hat{\alpha}_1 \rangle\Phi_1 +
\langle\hat{\mu}_0,\hat{\alpha}_2 \rangle \Phi_2\\
\Phi_1=\ln (1+|z_1|^2)-\ln (1+|z_2|^2), \quad \Phi_2= \ln
(1+|z_1|^2) + \ln (1+|z_2|^2).
\end{gather*}
Here the bilinear form on $\mathfrak{so}(4)$ is defined by $\langle
A,B \rangle = \frac{1}{2}\Tr AB$.
\end{ex}

\begin{prop}
The K\"{a}hlerian potential on each coadjoint orbit
$\mathcal{O}_{\mu_0}$ of a compact classical Lie group $\mathrm{G}$
has the following form
$$\Phi = \sum_k \langle {\mu}_0,{\alpha}_k \rangle \Phi_k, \qquad \Phi_k= a_{\alpha_k},$$ where
$\alpha_k$ is a simple root of $\mathfrak{g}$, and $a_{\alpha_k}$ is
the canonical coordinate corresponding to $H_{\alpha_k}\,{\in}\,
\mathfrak{h}$, and $\langle\cdot, \cdot \rangle$ denotes a bilinear
form on the dual space to~$\mathfrak{g}$.
\end{prop}

\begin{rem}
If $\mu_0$ satisfies the \textbf{integer condition}
$$2\frac{\langle {\mu}_0,{\alpha}_k \rangle }{\langle {\alpha}_k,{\alpha}_k\rangle} \in
\Integer,$$ then the orbit through $\mu_0$ can be quantized. In
other words, there exists an irreducible unitary representation of
$\mathrm{G}$
 in the space of holomorphic sections on the
orbit. Each section serves as a quantum state.
\end{rem}

%==============================================================
\section{Cohomology Rings of Coadjoint Orbits}\label{s:topology}
In the last section we examine the cohomology rings of coadjoint
orbits of compact semisimple Lie groups. A.~Borel \cite{borel53}
proved that all forms of odd degrees on the orbit are precise.
Therefore, we are interested in the forms of even degrees. In order
to introduce a basis for the cohomology ring it is sufficient to
find a basis for the cohomology group $H^2$.

In the case of a generic coadjoint orbit of a compact semisimple Lie
group $\mathrm{G}$, the following formula is available
$$b^0 + b^2 + \dots + b^{2n} = \ord \mathrm{W}(\mathrm{G}),$$ where $b^k$ denotes the
Betti number of a cohomology group $H^k$. In the case of a
degenerate orbit, one has to modify the formula as
$$b^0 + b^2 + \dots + b^{2m} = \frac{\ord
\mathrm{W}(\mathrm{G})}{\ord \mathrm{W}(\mathrm{G}_{\mu_0})},$$
where $\mathrm{G}_{\mu_0}$ is the stability subgroup at $\mu_0$.

\begin{ex}
In the case of group $\mathrm{SU}(2)$, we have the only  type of
orbits: $\mathcal{O}^{\mathrm{SU}(2)}$ of dimension 2. The Weyl
group $\mathrm{W}(\mathrm{SU}(2))$ also has dimension 2. Therefore,
the cohomology ring consists of two cohomology groups, each of
dimension~1:
$$ H^{\ast} = H^0 \oplus H^2, \qquad  1+1=2. $$

In the case of group $\mathrm{SU}(3)$, we have two types of orbits:
a generic one $\mathcal{O}^{\mathrm{SU}(3)}$ of dimension 6, and a
degenerate one $\mathcal{O}_d^{\mathrm{SU}(3)}$ of dimension 4. In
the case of a generic orbit, the Weyl group has dimension 6, and the
cohomology ring is
$$H^{\ast} = H^0 \oplus H^2 \oplus H^4 \oplus H^6,\qquad 1+2+2+1=6.$$

For a degenerate orbit we have $\frac{\ord
\mathrm{W}(\mathrm{G})}{\ord \mathrm{W}(\mathrm{G}_{\mu_0})} =3$,
and the cohomology ring is
$$H^{\ast} = H^0 \oplus H^2 \oplus H^4,\qquad 1+1+1=3.$$
\end{ex}

Recall the well-known Leray-Hirsch theorem.
\begin{theorem*}[Leray-Hirsch]
Suppose $\mathcal{E}$ is a fibre bundle over a base $\mathcal{M}$
with a fibre~$\mathcal{F}$, and $\omega_1$, $\omega_2$, \ldots
$\omega_r$ are cohomology classes on $\mathcal{E}$ that being
restricted to each fibre give its cohomologies. Then
$$H^{\ast}(\mathcal{E}) = H^{\ast}(\mathcal{M}) \otimes
H^{\ast}(\mathcal{F}).$$
\end{theorem*}
Apply the theorem to an orbit $\mathcal{O}$ regarded as a fibre
bundle over an orbit $\mathcal{O}_1$ with an orbit $\mathcal{O}_2$
as a fibre, that is $\mathcal{O} = \mathcal{E}(\mathcal{O}_1,
\mathcal{O}_2,\pi)$. The cohomology ring of $\mathcal{O}$ is a
tensor product of the cohomology rings of the base and the fiber:
$$H^{\ast}(\mathcal{O}) = H^{\ast}(\mathcal{O}_1) \otimes
H^{\ast}(\mathcal{O}_2).$$

Conversely, if one finds coherent cohomology classes on
$\mathcal{O}_1$ and $\mathcal{O}_2$, then one can construct the
cohomology ring of $\mathcal{O}$ by the latter formula. It means,
the cohomology ring of a generic orbit can be deriven from the
cohomology rings of a degenerate orbit and a generic orbit of a
group of less dimension.

\begin{ex}
We continue to deal with the group $\mathrm{SU}(3)$. It was shown
that
$$\mathcal{O}^{\mathrm{SU}(3)} = \mathcal{E}(\mathcal{O}_d^{\mathrm{SU}(3)}
, \mathcal{O}^{\mathrm{SU}(2)},\pi).$$ Then the cohomology ring of
$\mathcal{O}^{\mathrm{SU}(3)}$ is the tensor product of the
cohomology rings of the orbits $\mathcal{O}_d^{\mathrm{SU}(3)}$ and
$\mathcal{O}^{\mathrm{SU}(2)}$:
\begin{multline*}
H^{\ast}(\mathcal{O}^{\mathrm{SU}(3)}) =
(H^0 \oplus H^2 \oplus H^4) \otimes (H^0 \oplus H^2) = \\
= H^0 \otimes H^0 \oplus \underbrace{H^0 \otimes H^2 \oplus H^2
\otimes H^0}_{H^2(\mathcal{O}^{\mathrm{SU}(3)})} \oplus
\underbrace{H^2 \otimes H^2 \oplus H^4 \otimes
H^0}_{H^4(\mathcal{O}^{\mathrm{SU}(3)})} \oplus H^4 \otimes H^2.
\end{multline*}
Obviously, the cohomology groups $H^2$ and $H^4$ of
$\mathcal{O}^{\mathrm{SU}(3)}$ both have dimension 2. Moreover, from
the previous expression we can see the structure of a basis for
$H^2$:
$$H^2(\mathcal{O}^{\mathrm{SU}(3)}) = H^0(1)\otimes H^2(2) \oplus
H^2(1)\otimes H^0(2),$$  where  1 denotes
$\mathcal{O}^{\mathrm{SU}(3)}_{d}\simeq \CP^2$, and 2 denotes
$\mathcal{O}^{\mathrm{SU}(2)}\simeq \CP^1$.
\end{ex}

At the same time, a suitable basis for $H^2$ can be obtained from
K\"{a}hlerian potentials on coadjoint orbits of a group. As shown in
the previous section, all two-forms on the orbits of a compact
classical Lie group $\mathrm{G}$ have the form
$$\omega = \sum_k \i c_k \sum_{\alpha,\beta} \frac{\partial^2
\Phi_k}{\partial z_{\alpha}
\partial \bar{z}_{\beta}} \, dz_{\alpha} \wedge d\bar{z}_{\beta},
\qquad k=1,\, \dots, \, \dim\mathrm{T},$$ where $\Phi_k$ coincides
with the canonical coordinate $a_{\alpha_k}$ corresponding to
$H_{\alpha_k}\,{\in}\, \mathfrak{h}$. Obviously, $\dim H^2 = \dim
\mathrm{T}=l$. Consequently, one can find precisely $l$ two-forms
that give a basis for $H^2$.

The standard way to generate a basis for $H^2$ is the following. Let
$H_2$ be the homology group adjoint to $H^2$. By $[\gamma]$ we
denote a class of two-cycles, which can be represented as spheres.
The sphere is an orbit of a subgroup $\mathrm{SU}_{\alpha}(2)$:
$$\mathrm{SU}_{\alpha}(2) \simeq \exp\{H_{\alpha},\, (X_{\alpha}-X_{-\alpha}),\,
\i (X_{\alpha}+X_{-\alpha})\},\quad \alpha\in \Delta_{+}.$$ Suppose
we find $l$ independent two-cycles connected with the simple roots
of $\mathfrak{g}$, we denote them by $\gamma_i$. The basis for $H^2$
consists of two-forms $\omega_j$ such that
\begin{equation}\label{FormCond}
\int_{\gamma_i} \omega_j = \delta_{ij},
\end{equation}
where $\delta_{ij}$ is the Croneker symbol.

\begin{ex}
We consider coadjoint orbits of $\mathrm{SU}(3)$ as an example. Let
simple roots of $\mathfrak{su}(3)$ be as follows:
$\alpha_1=\diag(\i,-\i,0)$ and $\alpha_2=\diag(0,\i,-\i)$. Then
independent two-cycles are generated by the following dressing
matrices $$\hat{u}_1 = \begin{pmatrix} \frac{1}{\sqrt{1+|z_1|^2}} &
\frac{-\bar{z}_1}{\sqrt{1+|z_1|^2}} & 0 \\
\frac{z_1}{\sqrt{1+|z_1|^2}} & \frac{1}{\sqrt{1+|z_1|^2}} & 0 \\
0&0&1 \end{pmatrix},\qquad \hat{u}_2 = \begin{pmatrix} 1&0&0\\ 0&
\frac{1}{\sqrt{1+|z_2|^2}} &
\frac{-\bar{z}_2}{\sqrt{1+|z_2|^2}} \\
0 & \frac{z_2}{\sqrt{1+|z_2|^2}} & \frac{1}{\sqrt{1+|z_2|^2}}
\end{pmatrix},$$ which are obtained from the dressing matrix
$\hat{u}$ by assigning $z_2=z_3=0$ or $z_1=z_3=0$, respectively. The
 two-forms $\omega_j$ satisfying \eqref{FormCond} are
\begin{gather*}
\omega_j = \frac{1}{2\pi}\,\sum_{\alpha,\beta} \frac{\partial^2
\Phi_j}{\partial z_{\alpha} \partial \bar{z}_{\beta}} \, dz_{\alpha}
\wedge d\bar{z}_{\beta}, \quad j=1,\,2 \\ \Phi_1 = \ln
(1+|z_1|^2+|z_3-z_1 z_2|^2),\qquad \Phi_2 = \ln (1+|z_2|^2+|z_3|^2).
\end{gather*}
They form a basis for $H^2(\mathcal{O}^{\mathrm{SU}(3)})$.
\end{ex}

\section{Conclusion}
In this paper we develop a unified approach to solutions of the
announced problems for a coadjoint orbit of a compact semisimple
classical Lie group $\mathrm{G}$. The problems are the following: an
explicit parameterization of the orbit, obtaining a K\"{a}hlerian
structure, introducing basis forms for the cohomology group of the
orbit.  The key role belongs to the subgroup $\mathrm{A}$ in an
Iwasawa decomposition, this is the real abelian subgroup of a
complexification  of the group $\mathrm{G}$. The subgroup
$\mathrm{A}$ determines a K\"{a}hlerian potential on each orbit and
a suitable basis for the cohomology group $H^2$ of the orbit.

Our investigation concerns classical (matrix) Lie groups. The same
problems in the general case remain of current importance.

\section*{Acknowledgements}

We would like to thank professor A.~Arvanitoyergos for fruitful
debate and proper references, professor G.~Vilasi for valuable
discussion, professors I.~Mladenov and M.~Hadzhilazova for their
hospitality.

This work was partly supported by grant DFFD.Ukr F16/457-2007

\end{document}